\documentclass[secnum,leqno]{rims-bessatsu}
\usepackage{amssymb, amsmath}
\usepackage{graphics}
\usepackage[dvips]{graphicx}
\usepackage{eepic}
\usepackage{epic}
\newtheorem{thm}{Theorem}[section]

\theoremstyle{definition}
\newtheorem{dfn}[thm]{Definition}
\newtheorem{exa}[thm]{Example}
\theoremstyle{remark}
\title{Stochastic differential equations 
related to \\ random matrix theory}
\TitleHead{SDEs related to RMT}          
\author{Hirofumi \textsc{Osada}\footnote{
Faculty of Mathematics, Kyushu University,
Fukuoka 819-0395 Japan.
\newline e-mail: \texttt{osada@math.kyushu.ac.jp}}
          ~and Hideki \textsc{Tanemura}\footnote{Department of Mathematics and Informatics,
Faculty of Science, Chiba University, 
1-33 Yayoi-cho, Inage-ku, Chiba 263-8522, Japan. 
\endgraf e-mail: \texttt{tanemura@math.s.chiba-u.ac.jp}}}
\AuthorHead{H. Osada and H. Tanemura}         
\classification{15B52, 30C15, 47D07, 60G55, 82C22}
\support{Supported by a Grant-in-Aid for Scientific Research (KIBAN-A, No.24244010) and a Grant-in-Aid for Scientific Research (KIBAN-C, No.15K04916)}  
\keywords{\textit{Random matrices, Interacting particle systems, Stochastic differential equations}:}         
\VolumeNo{x}           
\YearNo{201x}           
\PagesNo{000--000}      
\communication{Received April 20, 201x.
Revised September 11, 201x.}      
\def\C{\mathbb{C}}

\def\I{\mathbb{I}}

\def\N{\mathbb{N}}

\def\R{\mathbb{R}}
\def\W{\mathbb{W}}

\def\d{{\bf d}}

\def\1{{\bf 1}}

\def\X{{\bf X}}

\def\cE{{\cal E}}

\def\x{{\bf x}}
\def\X{{\bf X}}

\def\Ai{{\rm Ai}}

\def\mM{\mathfrak{M}}

\def\={\stackrel{\rm (law)}{=}}

\begin{document}
%

\maketitle

\begin{abstract}      
In this note we review recent results on existence and uniqueness of solutions of infinite-dimensional stochastic differential equations describing interacting Brownian motions on $\R^d$.

\end{abstract}

\section{Introduction}

Let $\X^N(t)=(X_j^N(t))_{j=1}^{N}$ be a solution of the stochastic differential equation (SDE) 
\begin{align}&
\label{:1}
&&dX_j^N(t)=dB_j(t)+\frac{\beta}{2}\sum_{k=1,k\not=j}^N \frac{dt}{X_j^N(t)-X_k^N(t)}
\end{align}
or the SDE with Ornstein-Uhlenbeck's type drifts 
\begin{align} \label{:2} &
 dX_j^N(t)=dB_j(t) - \frac{\beta}{4N} X_j^N(t)dt  
+\frac{\beta}{2}\sum_{k=1,k\not=j}^N \frac{dt}{X_j^N(t)-X_k^N(t)}
,\end{align}
where $B_j(t), j=1,2,\dots, N$ are independent one-dimensional Brownian motions. These are called Dyson's Brownian motion models with parameters $\beta>0$ \cite{Dys62}. 
They were introduced to understand the statistics of eigenvalues of random matrix ensembles as distributions of particle positions in one-dimensional Coulomb gas systems with log-potential.

The solution of \eqref{:2} is a natural reversible stochastic dynamics with respect to 
$  \check{\mu}_{{\sf bulk},\beta}^N $:
\begin{equation}
\label{:1.3}
\check{\mu}_{{\sf bulk},\beta}^N(d\x_N)=\frac{1}{Z}h_N(\x_N)^\beta 
e^{-\frac{\beta}{4N}|\x_N|^2}
d\x_N,
\end{equation}
where $d\x_N=dx_1dx_2\cdots dx_N$, $\x_N=(x_i)\in\mathbb{R}^N$, and 
$$\displaystyle{h_N(\x_N)=\prod_{i<j}^N|x_i-x_j|}
.$$  
Throughout, $Z$ denotes a normalizing constant. 
Gaussian ensembles are called Gaussian orthogonal/unitary/symplectic ensembles (GOE/GUE/GSE) according to their invariance under conjugation by orthogonal/unitary/symplectic groups, which correspond to the inverse temperatures $\beta=1,2$ and $4$, respectively \cite{Meh04, AGZ10}.
It is natural to believe that the $N$-limit of the process $\X^N(t)$ solves the infinite-dimensional stochastic differential equation (ISDE)
\begin{equation}
\label{sin}
dX_j(t)=dB_j(t)+\frac{\beta}{2}\lim_{r\to\infty}\sum_{\substack{k=1,k\not=j \\ |X_k(t)|<r} }^{\infty}\frac{dt}{X_j(t)-X_k(t)}.
\end{equation}
The result was not proved rigorously until a few years ago when it was shown for $\beta=2$ in \cite{o-t.sm}, for $\beta=1,2,4$ in \cite{k-o.sgap}, and for $\beta\ge 1$ in \cite{tsai} .

Set $Y_j^N(t)=N^{1/6}(X_j^N(t)-2\sqrt{N})$, $j=1,2,\dots.N$ for the solution $\X^N$ of (\ref{:2}).
It has also been shown that
the $N$-limit of the process $\mathbf{Y}^N(t)$ solves the ISDE
\begin{equation}
\label{Ai} 
dY_j(t)=dB_j(t)
+\frac{\beta}{2}\lim_{r\to\infty}\left\{\sum_{\substack{k=1, k\not=j \\ |Y_k(t)|<r}}^{\infty}\frac{1}{Y_j(t)-Y_k(t)}-\int_{-r}^r \frac{\widehat{\rho}(x)dx}{-x}\right\}dt,
\end{equation}
with $\widehat{\rho}(x)=\pi^{-1}\sqrt{-x}\mathbf{1}(x<0)$, for $\beta=2$ \cite{o-t.sm} and for $\beta=1,2,4$ \cite{k-o.sgap}.

One of the key parts of proving the above results is  the existence and uniqueness of solutions of an ISDE  of the form
\begin{equation}\label{ISDE}
dX_j(t)=dB_j(t)
-\frac{1}{2}\nabla\Phi (X_j(t))dt
-\frac{1}{2}\sum_{k=1,k\not=j}^{\infty}\nabla\Psi(X_j(t),X_k(t))dt
\end{equation}
with free potential $\Phi$ and interaction (pair) potential $\Psi$. In ISDEs (\ref{sin}) and (\ref{Ai}), $\Psi$ is given by the log pair potential $-\beta\log |x-y|$.
The present note is a short summary of results on existence and uniqueness of solutions for ISDE (\ref{ISDE}).

\section{Quasi-Gibbs measure}

Let $S$ be a closed set in $\R^d$ such that $0\in S$ and $\overline{S^{\rm int}}=S$, where $S^{\rm int}$ denotes the interior of $S$. 
The configuration space $\mM$ of unlabelled particles is given by
\begin{eqnarray}
\mM &=& \Big\{\xi: 
\mbox{$\xi$ is a nonnegative integer valued Radon measure in $S$}
\Big\}
\\
&=&\Big\{ \xi(\cdot)=\displaystyle{\sum_{j\in \I}} \delta_{x_j}(\cdot):
\sharp \{j\in \I : x_j \in K\} < \infty, \mbox{ for any $K$ compact} \Big\},
\nonumber
\end{eqnarray}
where $\I$ is a countable set and $\delta_a$ is the Dirac measure at $a\in S$.
Thus $\mM$ is a Polish space with the vague topology.
We also introduce a subset $\mM_{\rm s.i}$ of $\mM$:
\begin{equation}
\mM_{\rm s.i.}=\{\xi\in\mM : \xi(\{x\})\le 1 \ \mbox{for all $x\in  S$, $\xi(S)=\infty$}\},
\end{equation}
that is, the set of configurations of an infinite number of particles without collisions.
For Borel measurable functions $\Phi : S\to \R\cup\{\infty\}$ and $\Psi : S\times S\to \R\cup\{\infty\}$ 
and a given increasing sequence $\{b_r\}$ of $\N$, we introduce the Hamiltonian
\begin{equation}
H_r(\xi)=H_r^{\Phi,\Psi}(\xi)=\sum_{x_j\in S_r}\Phi(x_j)
+\sum_{x_j,x_k\in S_r, j<k}\Psi(x_j,x_k),
\quad \xi=\sum_{j\in\I} \delta_{x_j},
\end{equation}
where $S_r=\{x\in S : |x|<b_r\}$.
We call $\Phi$ a free potential, and call $\Psi$ an interaction potential.
Let $\Lambda_r^m$ be the restriction of a Poisson random measure with intensity measure $dx$ on $\mathfrak{M}_r^{m}=\{\xi\in\mM :\xi(S_r)=m\}$.
We define maps $\pi_r, \pi_r^c :\mM\to \mM$ such that $\pi_r(\xi)=\xi (\cdot \cap S_r)$ and $\pi_r^c(\xi)=\xi (\cdot \cap S_r^c)$.
For two measures $\nu_1, \nu_2$ on a measurable space $(\Omega,\mathcal{F})$ we write $\nu_1 \le \nu_2$ if $\nu_1(A) \le \nu_2(A)$ for any $A\in\mathcal{F}$. 
We can now state the definition of a {\it quasi-Gibbs measure} \cite{o.rm, o.rm2}.

\begin{dfn}
A probability measure $\mu$ on $\mM$ is said to be 
a $(\Phi,\Psi)$-quasi Gibbs measure if its regular conditional probabilities
$$
\mu_{r,\xi}^m (d\zeta)= \mu(d\zeta | \pi_r^c(\zeta)=\pi_r^c(\xi), \zeta(S_r)=m), \quad r,m\in\N,
$$
satisfy that, for $\mu$-a.s. $\xi$,
$$
c^{-1}e^{-H_r(\eta)}\Lambda_r^m(\pi_{S_r}\in d\eta) \le 
\mu^m_{r,\xi}(\pi_{S_r}\in d\eta)
\le c e^{-H_r(\eta)}\Lambda_r^m(\pi_{S_r}\in d\eta).
$$
Here, $c=c(r,m,\xi)$ is a positive constant depending on $r, m$, and $\xi$. 
\end{dfn}

It is readily seen that the quasi-Gibbs property is a generalized notion of the canonical Gibbs property.
If $\mu$ is a $(\Phi, \Psi)$-quasi Gibbs measure, then $\mu$ is also a $(\Phi+\Phi_{\rm loc.bdd}, \Psi)$-quasi Gibbs measure for any locally bounded measurable function $\Phi_{\rm loc.bdd}$. In this sense, the notion of ``quasi-Gibbs'' seems to be robust. 
Information about the free potential of $\mu$ is determined from its {\it logarithmic derivative} \cite{o.isde}.

A function $f$ on $\mM$ is called a polynomial function if
\begin{equation}
\label{pol}
f(\xi) =Q \left( \langle \phi_1,\xi\rangle, \langle \phi_2,\xi\rangle, \dots, \langle \phi_\ell,\xi\rangle\right) 
\end{equation}
with $\phi_k \in C_c^{\infty}(S)$ and 
a polynomial function $Q$ on $\R^\ell$, where 
$\langle \phi, \xi\rangle = \int_S\phi(x)\xi(dx)$ and 
$ C_c^{\infty}(S)$ is the set of smooth functions with compact support. 
We denote by $\mathcal{P} $ the set of all polynomial functions on $\mM$. 

\begin{dfn}
We call $\d^\mu \in L^1_{loc}(S \times\mM,\mu^{[1]})$ 
the logarithmic derivative of $\mu$ if  
$$
\int_{S\times\mM}\d^\mu(x,\eta)f(x,\eta)d\mu^{[1]}(x,\eta)=-\int_{S \times\mM}\nabla_xf(x,\eta)d\mu^{[1]}(x,\eta)
$$
is satisfied for 
$f\in C_c^\infty(S)\otimes{\cal P}$.
Here
$\mu^{[k]}$ is the Campbell measure of $\mu$
$$
\mu^{[k]}(A\times B)=\int_A\mu_{\x}(B)\rho^k(\mathbf{x})d\x,
\quad
A\in{\cal B}(S^k), B\in{\cal B}(\mathfrak{M}),
$$
$\mu_{\x}$ is the reduced Palm measure conditioned at
 $\mathbf{x}\in S^k$
\begin{equation}
\mu_{\x}=\mu\left(
\cdot - \sum_{j=1}^k \delta_{x_j} \Bigg{|} \xi(x_j)\ge 1 \mbox{ for $j=1,2,\dots, k $}
\right),
\end{equation}
and $\rho^k$ is the $k$-correlation function for $k\in\N$.
\end{dfn}

Quasi-Gibbs measures inherit the following property from canonical Gibbs measures \cite[Lemma 11.2]{ot.2}.
Let $\mathcal{T}(\mM)$ be the tail $\sigma$-field
$$
\mathcal{T}(\mM)= \bigcap_{r=1}^\infty \sigma(\pi_r^c)
$$
and let 
$\mu_{\rm Tail}^{\xi}$ be the regular conditional probability  defined as
\begin{equation}
\mu_{\rm Tail}^{\xi}=\mu(\cdot| \mathcal{T}(\mM))(\xi).
\end{equation}
Then the following decomposition holds:
\begin{equation}
\mu(\cdot)=\int_{\mM}\mu_{\rm Tail}^{\xi}(\cdot)\mu(d\xi).
\end{equation}
Furthermore, there exists a subset $\mM_0$ of $\mM$ satisfying $\mu(\mM_0)=1$ and, for all $\xi, \eta\in\mM_0$:
\begin{eqnarray}
&&\mu_{\rm Tail}^{\xi}(A)\in\{0,1\} \quad \mbox{ for all $A\in\mathcal{T}(\mM)$},\\
&&\mu_{\rm Tail}^{\xi}(\{\zeta\in\mM : \mu_{\rm Tail}^{\xi}=\mu_{\rm Tail}^{\zeta}\})=1,
\\
&&
\mbox{$\mu_{\rm Tail}^{\xi}$ and $\mu_{\rm Tail}^{\eta}$ are mutually singular on $\mathcal{T}(\mM)$ if $\mu_{\rm Tail}^{\xi}\not= \mu_{\rm Tail}^{\eta}$}.
\end{eqnarray}

\section{General theory of solutions of ISDEs}

A polynomial function $f$ on $\mM$ is a {\it local} function, that is, a function satisfying $f(\xi)=f(\pi_r(\xi))$ for some $r\in\N$.
When $\xi\in\mM_r^m$, $m\in\N\cup\{0\}$ and $\pi_r(\xi)$ is represented by $\sum_{j=1}^m\delta_{x_j}$, we can regard 
$f(\xi)=f(\sum_{j=1}^m\delta_{x_j})$ as a permutation invariant smooth function on $S_r^m$. 
For $f,g\in\mathcal{P}$, define 
$$
\mathbb{D}(f,g)(\xi)= \frac{1}{2}
\sum_{j=1}^\infty \nabla_{x_j}f(\xi)\cdot \nabla_{x_j}g(\xi)
.$$
For a probability $\mu$ on $\mM$,
we denote by $ L^2(\mM,\mu)$ the space of square integrable functions on $\mM$ with
the inner product $\langle \cdot, \cdot \rangle_\mu$ and the norm $\|\cdot\|_{L^2(\mM,\mu)}$.
We consider the bilinear form $(\cE^\mu, \mathcal{P}^\mu)$ on $L^2(\mM,\mu)$ defined by
\begin{eqnarray}
&&\mathcal{E}^\mu (f,g)= \int_{\mM}\mathbb{D}(f,g)d\mu,
\quad
\mathcal{P}^\mu=\{f\in \mathcal{P}: \|f\|_1^2 <\infty \}
\label{def:D}
,\end{eqnarray}
where $\|f\|_1^2 \equiv \cE^\mu (f,f)+\|f \|_{L^2(\mM,\mu)}^2$. 

We make the following assumptions

\vskip 1mm

\noindent (A.0) \quad $\mu$ has a locally bounded $n$-correlation function $\rho^n$ for each $n\in\N$.

\noindent (A.1) \quad
There exist upper semi-continuous functions $\Phi_0:S\to\mathbb{R}\cup\{\infty\}$ and $\Psi_0 :S\times S \to \mathbb{R}\cup\{\infty\}$ that are locally bounded from below, 
 and $c>0$ such that
$$
c^{-1}\Phi_0(x) \le \Phi (x) \le c\Phi_0(x),
\quad
c^{-1}\Psi_0(x,y) \le \Psi (x,y) \le c\Psi_0(x,y).
$$
\noindent (A.2) \quad There exists a $T>0$ such that for each $R>0$
$$
\liminf_{r\to\infty}
\mathrm{Erf}\left(\frac{r}{(r+R)T}\right)\int_{|x|\le r+R} \rho^1(x)dx =0,
$$
where 
$\mathrm{Erf}(t)=(2\pi)^{-1/2}\int_t^\infty e^{-x^2/2}dx$.

Note that $\mathcal{P}^{\mu}=\mathcal{P}$ and $(\mathcal{E}^\mu,\mathcal{P}^\mu)=(\mathcal{E},\mathcal{P})$ under condition (A.0).

\begin{thm}[{\cite{o.isde, o.rm, o.rm2,o.tp, o-t.core}}]\label{existence}
Suppose that $\mu$ is a $(\Phi,\Psi)$-quasi Gibbs measure satisfying (A.0) and (A.1).
Then

\noindent
(i) $(\mathcal{E},\mathcal{P})$ is closable and its closure $(\mathcal{E}^{\mu},\mathcal{D}^\mu)$ is a quasi regular Dirichlet form and there exists the diffusion process $(\Xi(t), P_{\mu}^{\xi})$ associated with $(\mathcal{E}^{\mu},\mathcal{D}^\mu)$.
\\
(ii) Furthermore, assume conditions (A.2) and (A.3): 
\\
\noindent (A.3)
$\mathrm{Cap}^\mu((\mM_{\rm s.i})^c)=0$ and $\mathrm{Cap}^\mu(\xi(\partial S)\ge 1)=0$, 
\\
where $\mathrm{Cap}^\mu$ is the capacity of the Dirichlet form.
If there exists a logarithmic derivative $\d^\mu$,
then there exists $\tilde{\mathfrak{M}}\subset \mathfrak{M}$ such that $\mu(\tilde{\mathfrak{M}})=1$,
and for any $\xi=\sum_{j\in\mathbb{N}}\delta_{x_j}\in\tilde{\mathfrak{M}}$, there  exists an $S^\mathbb{N}$-valued continuous process $\mathbf{X}(t)=(X_j(t))_{j=1}^\infty$ satisfying 
${\bf X}(0)=\x=(x_j)_{j=1}^\infty$ and
$$
dX_j(t)=dB_j(t)+\frac{1}{2}\d^\mu
\bigg(X_j(t),\sum_{k:k\not=j}\delta_{X_k(t)}\bigg) dt, \quad j\in\N.
$$
\end{thm}

Let $\mathfrak{l}$ be  a label map from $\mM_{\rm s.i.}$ to $S^{\N}$, that is, 
for each $\xi\in\mM_{\rm s.i.}$, $\mathfrak{l}(\xi)= (\mathfrak{l}(\xi)_j)_{j=1}^\infty\in S^{\N}$ satisfies
$\xi= \sum_{j=1}^{\infty}\delta_{\mathfrak{l}(\xi)_j}$.
The map $\mathfrak{l}$ can be lifted to the map from $C([0,\infty),\mM_{s.i.})$ to $C([0,\infty),S^{\N})$.
For $\Xi \in C([0,\infty),\mM_{\rm s.i.})$ we put
$$
\Xi^{\diamond m}(t)=\sum_{j=m+1}^{\infty}\delta_{X_j(t)}
$$ for each $m\in\N$, where $(X_j)_{j=1}^{\infty}=\mathfrak{l}(\Xi)\in C([0,\infty),S^{\N})$.
We make the following assumption.

\vskip 3mm

\noindent (A4) 
There exists a subset $\mM_{\rm SDE}$ of $\mM_{s.i.}$ such that
$$
P_{\mu}^{\xi}( \Xi(t)\in \mM_{\rm SDE})=1
\quad \mbox{for any $\xi\in\mM_{\rm SDE}$},
$$
and for each $\Xi \in C([0,\infty),\mM_{\rm SDE})$ and each $m\in\N$,
\begin{eqnarray}
dY_j^{(m)}(t)&=&dB_j(t)
-\frac{1}{2}\nabla\Phi (Y_j^{(m)}(t))dt
-\frac{1}{2}\sum_{k=1,k\not=j}^{m}\nabla\Psi(Y_j^{(m)}(t),Y_k^{(m)}(t))dt
\label{finite}\\
&&\qquad -\frac{1}{2} \int_{\mM}\nabla\Psi(Y_j^{(m)}(t), X(t)) \Xi^{\diamond m}(dX)dt, \quad 1\le j\le m,
\nonumber
\\
Y_j^{(m)}(0)&=& \mathfrak{l}(\Xi(0))_j, \quad 1\le j\le m,
\end{eqnarray}
has a unique strong solution $\mathbf{Y}^{(m)}=(Y_1^{(m)},Y_2^{(m)},\dots,Y_m^{(m)})$.

We also make the following assumptions about the probability measure $\mu$

\noindent (A5)
For each $r,T\in \N$, there exists a positive constant $c$ such that
$$
\int_S \mathrm{Erf}\left(\frac{|x|-r}{\sqrt{cT}}\right)\rho^1(x)dx <\infty.
$$

\vskip 3mm

\noindent (A6)
The tail $\sigma$-field $\mathcal{T}(\mM)$ is $\mu$-trivial, that is, $\mu(A)\in \{0,1\}$ for $A\in \mathcal{T}(\mM)$. 

\begin{dfn}
Let $\mu$ be a probability measure on $\mM$ and let $\Xi (t)$ be an $\mM$-valued process.
We say that $\Xi (t)$ satisfies the {\it $\mu$-absolute continuity condition} if 
$\mu\circ\Xi(t)^{-1}$ is absolutely continuous with respect to $\mu$ for $\forall t>0$.
We say that an $S^{\N}$-valued process $\X(t)$ satisfies the $\mu$-absolute continuity condition if
$\frak{u}(\X (t))$ satisfies the $\mu$-absolute continuity condition, where $\mathfrak{u}$ is the map from $S^{\N}$ to $\mM$ defined by $\mathfrak{u}((x_j)_{j=1}^\infty)=\sum_{j=1}^\infty \delta_{x_j}$.
\end{dfn}

Then we have the following theorem.

\begin{thm}[{\cite{ot.2}}]\label{strong}
Suppose that the assumptions in Theorem \ref{existence} are satisfied.
Furthermore assume (A4)--(A6). Then, for $\mu$-a.s. $\xi$, 
ISDE (\ref{ISDE}) with $\X(0)=\mathfrak{l}(\xi)$ has a strong solution satisfying the $\mu$-absolute continuity condition, and that pathwise uniqueness holds for ISDE (\ref{ISDE}) with the $\mu$-absolute continuity condition.
\end{thm}

\section{Applications}

Theorems \ref{existence} and \ref{strong} can be applied to quite general class of ISDEs.
In this section we give some important examples.

\begin{exa}
[Canonical Gibbs measures]
Let $S=\mathbb{R}^d$, $d\in\N$.
Assume that $\Phi=0$ and that $\Psi_0$ is a super stable and regular in the sense of Ruelle \cite{ruelle.2}, and is smooth outside the origin.
Let $\mu$ be a canonical Gibbs measure with the interaction $\Psi_0$.
Then its logarithmic derivative is
\begin{equation}
\d^\mu
\bigg(x,\sum_{k:k\not=j}\delta_{y_k}\bigg)
=-\sum_{k=1,k\not=j}^{\infty}\nabla \Psi_0 (x-y_k).
\end{equation}
Assume that (A.2) is satisfied. 
In the case $d\ge 2$, there exists a diffusion process associated with $\mu$ and the labeled process solves 
\begin{equation}
\label{Cann}
dX_j(t)=dB_j(t)
-\frac{1}{2}\sum_{k=1,k\not=j}^{\infty}\nabla\Psi_0(X_j(t)-X_k(t))dt.
\end{equation}
In the case $d=1$, $\Psi_0$ needs to be sufficient repulsive at the origin to satisfy (A.3). 

Assume that (A.5) is satisfied and that, for each $n\in\N$, there exist positive constants $c, c'$ satisfying
\begin{eqnarray}
&&\sum_{r=1}^\infty \frac{\int_{|x|>r}\rho^1(x)dx}{r^c} <\infty,
\\
&&\sum_{i,j=1}^d
\left|\frac{\partial^2}{\partial x_i \partial x_j} \Psi_0(x)\right| \le \frac{c'}{(1+|x|)^{c'+1}},
\end{eqnarray}
for all $|x|\ge 1/n$.
In \cite[Theorem 3.3]{ot.2} 
it was proved that, for $\mu$-a.s. $\xi$, ISDE (\ref{Cann}) with $\X(0)=\mathfrak{l}(\xi)$ has a strong solution satisfying the $\mu_{\rm Tail}^\xi$-absolute continuity condition, and that pathwise uniqueness holds for ISDE (\ref{ISDE}) with the $\mu_{\rm Tail}^\xi$-absolute continuity condition.
\end{exa}

\begin{exa}[Sine random point fields]
Let $\check{\mu}_{{\sf bulk},\beta}^N$ be the probability measure defined in (\ref{:1.3}).
We denote by $\mu_{{\sf bulk},\beta}^N$ the distribution of $\sum_{j=1}^N \delta_{x_j}$ under $\check{\mu}_{{\sf bulk},\beta}^N$.
For $\beta>0$ the existence of the limit of $\mu_{{\sf bulk},\beta}^N$ as $N\to \infty$ was shown in Val\'ko-Vir\'ag\cite{VV}. We denote the limit by $\mu_{{\sf bulk},\beta}$. 
In particular, when $\beta=2$, $\mu_{{\sf bulk},2}$ is the determinantal point process (DPP) with the sine kernel
\begin{equation}
\label{def:sin}
K_{\sin, 2}(x,y)=\frac{\sin ( x-y)}{\pi(x-y)} 
,\end{equation}
and when $\beta=1, 4$, it is a quaternion determinantal point process \cite{AGZ10}.
It was shown that $\mu_{{\sf bulk},\beta}$ for $\beta=1,2,4$ is a quasi-Gibbs measure in \cite{o.rm}, and that its logarithmic derivative is
\begin{equation}\label{log_bulk}
\d^\mu
\bigg(x,\sum_{k:k\not=j}\delta_{y_k}\bigg)
=\beta
\lim_{r\to\infty}\sum_{\substack{k:k\not=j \\ |y_k|<r}}
\frac{1}{x-y_k}
\end{equation}
in \cite{o.isde}.
In \cite[Theorem 3.1]{ot.2} it was shown that 
for $\mu_{{\sf bulk},\beta}$-a.s. $\xi$,
ISDE (\ref{sin}) with $\X(0)=\mathfrak{l}(\xi)$ has a strong solution satisfying the $\mu_{{\sf bulk},\beta,{\rm Tail}}^\xi$-absolute continuity condition, and that pathwise uniqueness holds for ISDE (\ref{sin}) with the $\mu_{{\sf bulk},\beta,{\rm Tail}}^\xi$-absolute continuity condition.
In the case $\beta=2$, the facts that $\mathcal{T}(\mM)$ is $\mu_{{\sf bulk},2}$-trivial and $\mu_{{\sf bulk},2,{\rm Tail}}^\xi =\mu_{{\sf bulk},2}$ were shown in \cite{o-o.tail}.

Tsai \cite{tsai} proved the existence and uniqueness of solutions of ISDE (\ref{sin}) for $\beta\ge 1$ by a different method. Thus it is conjectured that $\mu_{{\sf bulk},\beta}$ is a quasi-Gibbs measure and has a logarithmic derivative of the form (\ref{log_bulk}) for $\beta\ge 1$.
\end{exa}

\begin{exa}[Airy random point fields]
We denote by $\mu_{{\sf soft},\beta}^N$ the distribution of $\sum_{j=1}^N \delta_{N^{1/6}(x_j-2\sqrt{N})}$ under $\check{\mu}_{{\sf bulk},\beta}^N$.
For $\beta>0$, the existence of the limit of $\mu_{{\sf soft},\beta}^N$ as $N\to \infty$ was shown in Ram\'irez-Rider-Vir\'ag \cite{RRV}. We denote the limit by $\mu_{{\sf soft},\beta}$. 
In particular, when $\beta=2$, $\mu_{{\sf soft},2}$ is the DPP with the Airy kernel
\begin{equation}
\label{def:Ai}
K_{\Ai, 2}(x,y)=\frac{\Ai(x)\Ai'(y)-\Ai'(x)\Ai(y)}{x-y},
\end{equation}
where $\Ai$ denotes the Airy function and $\Ai'$ its derivative \cite{Meh04}.
When $\beta=1, 4$, it is a quaternion determinantal point process \cite{AGZ10}.
In the cases $\beta=1,2,4$, it has been proved that the random point field is 
quasi-Gibbsian \cite{o.rm2},
and that its logarithmic derivative is
\begin{equation}\label{log_soft}
\d^\mu
\bigg(x,\sum_{k:k\not=j}\delta_{y_k}\bigg)
=\beta\lim_{r\to\infty}\left\{\sum_{\substack{k:k\not=j \\ |y_k|<r}}\frac{1}{x-y_k}-\int_{-r}^r \frac{\widehat{\rho}(x)dx}{-x}\right\},
\end{equation}
and for $\mu_{{\sf soft},\beta}$-a.s. $\xi$,
ISDE (\ref{Ai}) with $\X(0)=\mathfrak{l}(\xi)$ has a strong solution satisfying the $\mu_{{\sf soft},\beta,{\rm Tail}}^\xi$-absolute continuity condition, and pathwise uniqueness holds for ISDE (\ref{Ai}) with the $\mu_{{\sf soft},\beta,{\rm Tail}}^\xi$-absolute continuity condition \cite[Theorem 2.3]{o-t.airy}.
In the case $\beta=2$ the facts that $\mathcal{T}(\mM)$ is $\mu_{{\sf soft},2}$-trivial and that $\mu_{{\sf soft},2,{\rm Tail}}^\xi =\mu_{{\sf soft},2}$ were shown in \cite{o-o.tail}.

Determining whether $\mu_{{\sf soft},\beta}$ has the quasi-Gibbs property for general $\beta$ and finding its logarithmic derivative is (\ref{log_soft}) are interesting and important problems.
\end{exa}

\begin{exa}[Bessel random point field]
Let $S=[0,\infty)$ and $1\le \alpha < \infty$.
Let $\mu_{{\sf hard},2}$ be the determinantal point process with Bessel kernel
\begin{equation}
K_{J_{\alpha}}(x,y)=\frac{J_{\alpha}(\sqrt{x})\sqrt{y}J_{\alpha}'(\sqrt{y})- \sqrt{x}J_{\alpha}'(\sqrt{x})J_{\alpha}(\sqrt{y})}{2(x-y)}.
\end{equation}
In \cite{h-o.bes} it was shown that $\mu_{{\sf hard},2}$ is a quasi-Gibbs measure and that the related process is the unique strong solution of the ISDE
$$
dX_j(t)=dB_j(t)+\left\{\frac{\alpha}{2X_j(t)}+\sum_{k=1,k\not=j}^\infty\frac{1}{X_j(t)-X_k(t)}\right\}dt
$$
with the $\mu_{{\sf hard},2}$-absolute continuity condition.
\end{exa}

\begin{exa}[Ginibre random point field]
Let $S=\R^2$ be identified as $\C$. Let $\mu_{{\sf Gin}}$ be the DPP
with the kernel $K_{\sf Gin} : \C\times \C\to\C$ defined by
\begin{equation}
K_{\sf Gin}(x,y)=\frac{1}{\pi}e^{-|x|^2/2-|y|^2/2}e^{x\overline{y}}.
\end{equation}
In \cite{o.rm} it was shown that $\mu_{{\sf Gin}}$ is a quasi-Gibbs measure, and in \cite{o.isde} that the related process is a solution of the ISDE
\begin{equation}\label{Gin}
dX_j(t)=dB_j(t)-X_j(t)dt+\lim_{r\to\infty}
\sum_{\substack{k:k\not=j \\ |X_k(t)|<r}}
\frac{X_j(t)-X_k(t)}{|X_j(t)-X_k(t)|^2}dt.
\end{equation}
The pathwise uniqueness of solutions of (\ref{Gin}) with the $\mu_{{\sf Gin}}$-absolute continuity condition was shown in \cite{ot.2}.
\end{exa}

\section{Remarks}

In the previous section we gave some examples of DPPs that are not canonical Gibbs measures but quasi-Gibbs measures. 
It is expected that quite general DPPs have the quasi-Gibbs property. We thus present examples of DPPs related to random matrix theory or non-colliding Brownian motions, whose quasi-Gibbs property have not been shown.

\begin{exa}[Pearcey process]
Consider $2N$ noncolliding Brownian motions, in which all particles start from the origin and $N$ particles end at $\sqrt{N}$ at time $t=1$,
and the other $N$ particles end at $-\sqrt{N}$ at $t=1$.
We denote the system by $(X_1^N(t), \dots, X_{2N}^N(t))$, $0\le t\le 1$.
When $N$ is very large, there is a cusp at $x_0^N=0$ when $t_0=\frac{1}{2}$,
that is, before time $t_0$ particles are in one interval with high probability,
while after time $t_0$ they are separated into two intervals by the origin.
We denote the distribution 
$$
\sum_{j=1}^{2N}\delta_{2^{3/2} (2N)^{1/4}X_j^N(\frac{1}{2})}
$$
on $\mM$ by $\mu_{{\sf pearcey}}^{N}$.
It was proved in Adler-Orantin-von Moerbeke \cite{AOM} that 
$$
\mu_{{\sf pearcey}}^{N} \to \mu_{{\sf pearcey}},
\quad \mbox{weakly as $N\to\infty$}
$$
and that
$\mu_{{\sf pearcey}}$ is the DPP $K_{\sf pearcey}(x,y)$ given by
$$
K_{\sf pearcey}(x,y)=\frac{P(x)Q''(y)-P'(x)Q'(y)+P''(x)Q(y)}{x-y}, \quad x,y\in\R,
$$
with
$$
Q(y)=\frac{i}{2\pi}\int_{-i\infty}^{i\infty} e^{-u^4/4 -uy}du
\quad\mbox{ and }\quad
P(x)=\frac{1}{2\pi i}\int_C e^{v^4/4+vx}dv,
$$
where the contour $C$ is given by the ingoing rays from $\pm\infty e^{i\pi/4}$ to $0$
and the outgoing rays from $0$ to $\pm\infty e^{-i\pi/4}$.
These integrals are known as Pearcey's integrals \cite{pea}.
\end{exa}

\vskip 3mm

\begin{exa}[Tacnode process]
Consider two groups of non-colliding pinned Brownian motions $(X_1^N(t), \dots, X_{2N}^N(t))$ in the time interval $0\le t \le 1$, where one group of $N$ particles starts and ends at $\sqrt{N}$ and the other group of $N$ particles starts and ends at $-\sqrt{N}$.
The distribution $( N^{1/6}X_1^N(\frac{1}{2}), N^{1/6}X_2^N(\frac{1}{2}),\dots, N^{1/6}X_{2N}^N(\frac{1}{2}))$
 on the Weyl chamber of type $A_{2N-1}$
$$
\W_{2N}=\Big\{\x=(x_1, x_2, \cdots, x_{2N}):
x_1 < x_2 < \cdots < x_{2N} \Big\},
$$
is given by
\begin{equation}
m_{{\sf tac}}^{2N}(d\x_{2N})=\frac{1}{Z}
\left[\det_{1\le i,j \le 2N}\left(
e^{ -2|x_i-a_j|^2}
\right)\right]^2,
\nonumber
\end{equation}
where $a_j=-\sqrt{N}$ for $1\le j \le N$ and $a_{j}=\sqrt{N}$ for $N+1\le j \le 2N$.
We denote the distribution of $\sum_{j=1}^{2N}\delta_{N^{1/6}x_j}$ under $m_{{\sf tac}}^{2N}$ by $\mu_{{\sf tac}}^{N}$.
It was proved in Delvaux-Kuijlaars-Zhang \cite{DKZ11}
and Johansson \cite{J12} that 
$$
\mu_{{\sf tac}}^{N} \to \mu_{{\sf tac}},
\quad \mbox{weakly as $N\to\infty$}
$$
and that
$\mu_{{\sf tac}}$ is the DPP with the correlation kernel
\begin{equation}
K_{{\sf tac}}(x,y)\equiv L_{{\sf tac}}(x,y)+ L_{{\sf tac}}(-x,-y),
\quad x,y\in\R,
\nonumber
\end{equation}
where
\begin{eqnarray}
&&L_{{\sf tac}}(x,y)= K_{\Ai, 2}(x,y)
\nonumber\\
&&\qquad+2^{1/3}\int_{(0,\infty)^2} dudv \ {\Ai}(y+2^{1/3}u)R(u,v){\Ai}(x+2^{1/3}v)
\nonumber\\
&&\qquad -2^{1/3}\int_{(0,\infty)^2} dudv \ {\Ai}(-y+2^{1/3}u){\Ai}(u+v){\Ai}(x+2^{1/3}v)
\nonumber\\
&&-2^{1/3}\int_{(0,\infty)^3} dudvdw \ {\Ai}(-y+2^{1/3}u)R(u,v){\Ai}(v+w){\Ai}(x+2^{1/3}w).
\nonumber
\end{eqnarray}
Here, $R(x,y)$ is the resolvent operator for the restriction of the Airy kernel to $[0,\infty)$, that is, the kernel of the operator
\begin{align}
&R= (I-K_{\Ai})^{-1}K_{\Ai}
\end{align}
on $L^2[0,\infty)$.

In \cite{DKZ11, J12} it was also shown that
$$
\Xi^N(t)\equiv\sum_{j=1}^{2N}\delta_{ N^{1/6}X_j(\frac{1}{2}+N^{-1/3}t)} \to \Xi(t), \quad \mbox{ as $N\to\infty$},
$$
in the sense of finite-dimensional distributions,
where $\Xi(t)$ is a reversible process with reversible measure $\mu_{{\sf tac}}$.
We expect that $\Xi(t)$ is the diffusion process associated with the Dirichlet form $(\mathcal{E}^{\mu_{{\sf tac}}}, \mathcal{D}^{\mu_{{\sf tac}}})$.
\end{exa}

\vskip 3mm


\end{document}